\documentclass[reqno]{amsart}
\begin{document}
\def\N{{\mathbb N}}
\def\Z{{\mathbb Z}}
\def\R{{\rm I\!R}}
\def\Re{{\rm Re}}
\def\GL{{\bf GL}}
\def\phi{\varphi}
\def\const{{\rm const}}

\def\L{{\mathcal L}}
\def\P{{\mathbb P}}
\def\Q{{\mathcal Q}}
\def\V{{\mathcal V}}
\def\A{\ov A}
\def\Chi{K}
\def\Zeta{Z}
\def\c{{\mathcal C}}
\def\x{{\bar x}}
\def\y{{\bar y}}
\def\z{{\bar z}}
\def\w{{\bar w}}
\def\f{{\bar f}}
\def\g{{\bar g}}
\def\t{{\bar t}}

\def\C{{\mathchoice {\setbox0=\hbox{$\displaystyle\rm C$}\hbox{\hbox to0pt{\kern0.4\wd0\vrule height0.9\ht0\hss}\box0}}{\setbox0=\hbox{$\textstyle\rm C$}\hbox{\hbox to0pt{\kern0.4\wd0\vrule height0.9\ht0\hss}\box0}} {\setbox0=\hbox{$\scriptstyle\rm C$}\hbox{\hbox to0pt{\kern0.4\wd0\vrule height0.9\ht0\hss}\box0}} {\setbox0=\hbox{$\scriptscriptstyle\rm C$}\hbox{\hbox to0pt{\kern0.4\wd0\vrule height0.9\ht0\hss}\box0}}}}

\def\codim{{\rm codim}}
\def\M{{\mathcal M}}
\def\ov{\overline}
\def\m{\mapsto}
\def\rk{{\rm rank\,}}
\def\id{{\sf id}}
\def\Aut{{\sf Aut}}
\def\CR{{\rm CR\ }}
\def\crd{\dim_{{\rm CR}}}
\def\crc{{\rm codim_{CR}}}
\def\reg{{\rm reg}}
\def\sing{{\rm Sing}}
\def\eps{\varepsilon}
\def\f{{\bar f}}
\def\w{{\bar w}}
\def\x{{\bar x}}
\def\cn{{\C^n}}
\def\cnn{{\C^{n'}}}
\def\ocn{\ov{\C^n}}
\def\ocnn{\ov{\C^{n'}}}
\emergencystretch9pt
\frenchspacing

\newtheorem{Th}{Theorem}[section]
\newtheorem{Def}{Definition}[section]
\newtheorem{Cor}{Corollary}[section]
\newtheorem{Satz}{Satz}[section]
\newtheorem{Prop}{Proposition}[section]
\newtheorem{Lemma}{Lemma}[section]
\newtheorem{Rem}{Remark}[section]
\newtheorem{Example}{Example}[section]

\title[Algebraicity of local holomorphisms]{Algebraicity of local holomorphisms between real-algebraic submanifolds of complex spaces.}
\author{Dmitri~Zaitsev}
\address{Mathematisches Institut, Universit\"at T\"ubingen,
        72076 T\"ubingen, Germany,
        E-mail address: dmitri.zaitsev@uni-tuebingen.de}
\begin{abstract}
We prove that a germ of a holomorphic map $f$ between $\C^n$ and $\C^{n'}$
sending one real-algebraic submanifold $M\subset\C^n$ into another
$M'\subset\C^{n'}$ is algebraic provided $M'$ contains no complex-analytic discs
and $M$ is generic and minimal. We also propose an algorithm 
for finding complex-analytic discs in a real submanifold. 
\end{abstract}
\maketitle

\section{Introduction}

We say that a real submanifold $M\subset\C^n\cong\R^{2n}$
is {\em algebraic} if there exist polynomials
$p_1,\ldots,p_d\in\C[z,\z]$ and an open subset $U\subset\C^n$ such that
\begin{equation}\label{Meq}
M= \{ z\in U : p_1(z,\z)=\cdots = p_d(z,\z)=0 \}.
\end{equation}
Similarly we say that a holomorphic map $f$ from a domain $D$ into $\C^{n'}$
is {\em algebraic} if there exist polynomials
$P_1,\ldots,P_l\in\C[z,z']$ and an open subset
$V\times V'\subset D\times\C^{n'}$ such that
$$\{(z,f(z)) : z\in V \} =
\{ (z,z')\in V\times V' : P_1(z,z')=\cdots=P_l(z,z')=0 \}.$$
Finally we say that $f$ is a {\em local holomorphism} (resp. local biholomorphism) from
a real submanifold $M\subset\C^n$ into another one $M'\subset\C^{n'}$,
if $f$ is a holomorphic map from a domain $D\subset\C^n$ into $\C^{n'}$
(resp. a biholomorphic map between domains $D\subset\C^n$ and $D'\subset\C^{n'}$)
such that $M\cap D\ne\emptyset$ and $f(M\cap D)\subset M'$ (resp. $f(M\cap D)=M'\cap D'$).

Poincar\'e~\cite{Po} was one of the first who studied 
algebraicity properties of local biholomorphisms between hypersurfaces.
He proved that a local biholomorphism between pieces of $3$-spheres in $\C^2$ is a rational map.
This result was extended by Tanaka~\cite{Ta} to the case of higher dimensional spheres.
An important step in understanding this phenomenon was done by Webster~\cite{W} 
who proved the algebraicity of local biholomorphisms $f$
between Levi-nondegenerate algebraic hypersurfaces $M$ and $M'$.
The last main progress in this direction was done recently by
Baouendi, Ebenfelt and Rothschild~\cite{BER}, where they prove
the algebraicity of $f$ for $M$ (and $M'$) generic, minimal and essentially finite.
These conditions turn out to be optimal in most cases
(see \cite{BER}, Theorem~2).

A typical example of a biholomorphism is an extension
of a biholomorphic map between smooth domains $\Omega,\Omega'\subset \C^n$
to a boundary point of $\Omega$ if possible. The algebraicity of such a biholomorphism
implies in particular a holomorphic extension to a dense subset
of the boundary $\partial\Omega$. Similarly the extension phenomena
for proper holomorphic maps lead to the study of local holomorphisms
between hypersurfaces of different dimensions.
An important step here was done by Huang~\cite{Hua}
who proved the algebraicity of holomorphisms 
between strictly pseudoconvex hypersurfaces.
See also \cite{W2}, \cite{TH}, \cite{Tu1}, \cite{F2}, \cite{S1}, \cite{S3}, \cite{S4}, \cite{SS}
for algebraicity results in other situations.

In this paper we consider the general case of real-algebraic
submanifolds $M\subset\cn$ and $M'\subset\cnn$ of arbitrary dimensions and codimensions.
Our main result is the following:

\begin{Th}\label{main}
Let $M\subset\cn$ and $M'\subset\cnn$ be real-algebraic submanifolds
such that the following is satisfied:
\begin{enumerate}
\item Every local holomorphism from $M$ into $\R$ is constant;
\item Every local holomorphism from $\C$ into $M'$ is constant.
\end{enumerate}
Then every local holomorphism from $M$ into $M'$ is algebraic.
\end{Th}

If $M$ connected, the first condition in Theorem~\ref{main} holds if and only if
$M$ is generic ($T_xM+iT_xM=T_x\C^n$) and minimal (in the sense of Tumanov~\cite{Tu})
at some point $x\in M$. In particular, it holds if the Levi cone of $M$ at a generic point
has a nonempty interior 
(if $M$ is a hypersurface, this means that $M$ is not locally biholomorphic to $\C^{n-1}\times\R$).
On the other hand, if this condition does not hold
and $M$ is homogeneous (i.e. the polynomials in (\ref{Meq}) can be chosen homogeneous with respect
to some weights), there always exist non-algebraic self-biholomorphisms of $M$ \cite{BER}.

The second condition in Theorem~\ref{main} means that there is no complex-analytic discs inside $M'$.
This holds, in particular, for boundaries of real-analytic bounded domains
(see \cite{DF0}) or for strictly pseudoconvex hypersurfaces. Therefore,
the algebraicity result of \cite{Hua} is contained in Theorem~\ref{main}.
On the other hand, if this condition does not hold and $\phi\colon\Delta\to M'$
is a complex-analytic disc, there always exist non-algebraic holomorphisms 
of the form $\phi\circ\pi$, where $\pi$ is a holomorphisms between a neighborhood
of a point of $M$ and $\Delta$.

Both conditions can be tested effectively.
It follows from Theorem~3 and Corollary~2.2.2 in \cite{BER} that
the first condition is valid if and only if $\pi\colon \M^m_z\to X_m$
is generically submersive, where $\M^m$ is defined in (\ref{Mm}) 
and $\M^m_z:=\{\xi\in\M^m : \xi_0=z\}$.
For the second condition we give an algorithm in \S\ref{algo}.

Our method of proving Theorem~\ref{main} 
is based on the geometric reflection principle
which includes the study of the families of Segre varieties (see \cite{W}) and of Segre sets 
(see \cite{BER}) associated to a real-analytic submanifold of $\C^n$. 
Recall that the Serge varieties are defined by $Q_w:=\{z : p_1(z,\w)=\cdots = p_d(z,\w)=0\}$,
where $M$ is supposed to be given by (\ref{Meq}) such that the differentials of
real and imaginary parts of the $p_i$'s are linearly independent on $M$.

The first main step is done by Proposition~\ref{psi} which establishes
an algebraic dependence between the jets of a local holomorphism $f$ at different points.
This can be seen as a generalization of the Schwarz reflection principle,
where the dependence is given by the standard reflection.
See \cite{DW}, \cite{BER}, \cite{BER1}, \cite{BER2} and \cite{Z2} for
the jet reflections of this kind in several complex variables.
The crucial assumption for this reflection was the essential finiteness
which guarantees that the image $f(Q_w)$ determines $Q'_{f(w)}$ up to finitely many possibilities.
In fact, one has locally $f(Q_w) = Q'_{f(w)}$ in all these cases.

The main difficulty in the general case of Theorem~\ref{main}
is that the above finite determinacy is no more valid. For example,
if $f$ is immersive and $M\subset\C^n$ and $M'\subset\cnn$ are hypersurfaces,
$\dim f(Q_w)=n-1$ whereas $\dim Q'_{f(w)}=n'-1$.
This makes impossible a direct reflection of jets as above.
Our idea here is to establish an algebraic relation between the jets
of $f$ at {\em three} (instead of two) different points $w$, $z$ and $w_1$. 
The $k$-jets $j^k_wf$ and $j^k_zf$ parametrize certain algebraic families
(denoted by $A$ and $B$) such that $f(w_1)$ is algebraically determined up to finitely many possibilities
by the {\em intersection of these families} (see the proof of Proposition~\ref{psi}).
This property is guaranteed by condition~2 in Theorem~\ref{main}.

An immediate consequence of Proposition~\ref{psi} is the algebraicity of $f$ along the Segre varieties.
The second main step is the iteration of the above reflection
which leads to the algebraicity of $f$ along the Segre sets $Q^s_w$
(the Segre sets can be defined inductively by $Q^1_w:=Q_w$, $Q^{s+1}_w:=\cup\{Q_z : z\in Q^s_w\}$).
Condition~1 in Theorem~\ref{main} guarantees that $\dim Q^s_w=n$ for some $s$.
Hence the algebraicity of $f$ along $Q^s_w$ is equivalent to its algebraicity.

F.~Meylan und A.~Sukhov informed the author about writing a manuscript jointly with B.~Coupet,
where they prove the statement of Theorem~\ref{main} in the case $M\subset\C^n$ is a hypersurface
and generalize it in other directions.

\section{Holomorphic and algebraic families}\label{fam}

All manifolds in the sequel are supposed to be complex unless otherwise stated.
All submanifolds and analytic subsets are complex and locally closed.
We say that a submanifold $\M$ of an algebraic variety $X$ is algebraic, if
for every $x\in\M$ there exist an affine neighborhood $U$ of $x$ and 
regular functions $P_1,\ldots,P_l$ in $U$ such that
$$ X\cap U = \{ z\in U: P_1(z)=\cdots=P_l(z)=0 \}.$$

For our construction we need maps which are
not necessarily defined everywhere.
A {\em partial (algebraic)} map $f\colon U\to V$
between (algebraic) manifolds is an (algebraic) map $f\colon S\to V$,
where $S\subset U$ is an (algebraic) submanifold of $U$.
When writing $f(u)$, we always assume $u\in S$.

\begin{Def}\label{family}
Let $U$ and $V$ be connected manifolds.
We call $F\subset U\times V$ a {\em holomorphic family}, if
there exists a splitting $V=V_1\times V_2$ such that
$$F=\{ (u,v_1,v_2)\in U\times V_1\times V_2 : v_2=\phi(u,v_1) \}$$
for some holomorphic map $\phi\colon U\times V_1\to V_2$.
If in addition $U$, $V$ and $F$ are algebraic
we call $F$ an {\em algebraic family}.
\end{Def}

To be more precise we have to speak of a family $(U,V,F)$.
However we shall identify $F$ with the triple $(U,V,F)$ for simplicity.
We write
$$F_u:=\{v\in V : (u,v)\in F\}\subset V, \quad u\in U,$$
for the fibers of $F$. Clearly a family is uniquely defined by its fibers.
For this reason we sometimes construct the fibers and speak
of the family $F_u\subset V$, $u\in U$.

\begin{Def}
Let $U$ and $V$ be manifolds and $F\subset U\times V$ be a subset.
A holomorphic family $\tilde F\subset \tilde U\times \tilde V$ is called
a {\em subfamily} of $F$, if $\tilde U$ and $\tilde V$
are open subsets of $U$ and $V$ respectively and
$\tilde F= F\cap (\tilde U\times \tilde V)$.
\end{Def}

In the sequel we shall frequently replace $F$ with a subfamily
satisfying special properties. If no ambiguity appears,
we denote the subfamily also by $F$ and speak
simply of ``replacing $F$ with a subfamily''.
The following is a direct consequence of the implicit function theorem:

\begin{Lemma}\label{submer}
Let $F\subset U\times V$ be a submanifold such that
the natural projection $\pi_{F,U}\colon F\to U$ is submersive at some point $x\in F$.
Then we can replace $F$ with a subfamily containing $x$.
\end{Lemma}

Recall that a real submanifold $M$ of a complex manifold $X$
is {\em generic}, if for all $x\in M$, $T_xM+iT_xM=T_xX$.
A generic submanifold cannot be included in a proper analytic subset of $X$.

\begin{Lemma}\label{image}
Let $F\subset U\times V$ be a holomorphic family,
$M\subset F$ a generic real submanifold,
$V'$ a complex manifold and $f\colon V\to V'$ a holomorphic map.
Then we can replace $F$ with a subfamily satisfying $M\cap F\ne\emptyset$ and
replace $V'$ with an open subset such that $f(V)\subset V'$ and
$f(F_u)\subset V'$, $u\in U$, is a holomorphic family.
\end{Lemma}

\begin{proof}
Set $g:=\id_U\times f\colon U\times V\to U\times V'$.
Since $M\subset F$ is generic, it cannot be included in the
analytic subset, where the rank of $g$ degenerates.
Then we can replace $F$ with a submanifold with a nonempty intersection with $M$
and such that $F':=g(F)$ is a submanifold of $V'$. Fix $x\in M$.
Then the surjectivity of $d_x\pi_{F,U}$ implies the surjectivity of $d_{g(x)}\pi_{F',U}$.
By Lemma~\ref{submer}, we can replace $F'$ with a subfamily containing $g(x)$.
It remains to replace $V$ with $f^{-1}(V')$.
\end{proof}

For $U$, $U'$ and $V$ arbitrary manifolds and
$F\subset U\times V$, $F'\subset U'\times V$
closed analytic subsets define
$$A_u(F,F'):= \{u'\in U': F'_{u'}\supset F_u\}, \quad u\in U,$$
and by $A(F,F')\subset U\times U'$ the corresponding family.

\begin{Lemma}\label{afg}
Let $F\subset U\times V$ be a holomorphic (resp. algebraic) family
and $F'\subset U'\times V$ be an analytic (resp. algebraic) subset.
Then $A(F,F')\subset U\times U'$
is a closed analytic (resp. algebraic) subset.
\end{Lemma}

\begin{proof}
Let $\phi\colon U\times V_1\to V_2$ be as in Definition~\ref{family}.
Since $F'$ is analytic (resp. algebraic),
it is locally defined by the vanishing
of holomorphic (resp. algebraic) functions
$f_1,\ldots,f_s$. Then $(u,u')\in A(F,F')$ is equivalent to the vanishing
of $f_i(u,v_1,\phi(u,v_1))$ for all $v_1\in V_1$ which is an analytic
(resp. algebraic) condition.
\end{proof}

Further we define for $A\subset U\times U'$ and
$G\subset U'\times V$:
$$B_u(A,G):= \bigcap_{u'\in A_u} G_{u'}\subset V , \quad u\in U,$$
and by $B(A,G)\subset U\times V$ the corresponding family
(in the case $A_u = \emptyset$ we set $B_u(A,G):= V$).

\begin{Lemma}\label{bag}
Let $A\subset U\times U'$ be a holomorphic (resp. algebraic) family
and $G\subset U'\times V$ an analytic (resp. algebraic) subset.
Then $B(A,G)\subset U\times V$ is a closed analytic (resp. algebraic) subset.
\end{Lemma}

\begin{proof}
By Definition~\ref{family} there exists a splitting $U'=U'_1\times U'_2$
and a holomorphic (resp. algebraic) map
$\phi\colon U\times U'_1\to U'_2$
such that
$$A=\{ (u,u'_1,u'_2)\in U\times U'_1\times U'_2 : u'_2=\phi(u,u'_1) \}.$$
Since $G$ is analytic (resp. algebraic),
it is locally defined by the vanishing
of holomorphic (resp. algebraic) functions $f_1,\ldots,f_s$.
Then $(u,v)\in B(A,G)$ is equivalent to the vanishing
of $f_i(u'_1,\phi(u,u'_1),v)$ for all $i$ and $u'_1\in U'_1$,
which is an analytic (resp. algebraic) condition.
\end{proof}

\section{Operations with jets}

Let $X$, $X'$ be arbitrary manifolds and $x\in X$ be an arbitrary point.
Recall that a $k$-jet at $x$ of a holomorphic map from $X$ into $X'$
is an equivalence class of germs at $x$ of local holomorphic maps between $X$ and $X'$
with fixed partial derivatives at $x$ up to order $k$.
Denote by $J^k_x(X,X')$ the space of all $k$-jets.
The union $J^k(X,X'):=\bigcup_{x\in X}J^k_x(X,X')$ carries
a natural fiber bundle structure over $X$.
For $f$ a holomorphic map from a neighborhood of $x$
into $X'$, we denote by $j^k_xf\in J^k_x(X,X')$ the corresponding $k$-jet.
If $X$ and $X'$ are algebraic manifolds,
$J^k_x(X,X')$ and $J^k(X,X')$ are also algebraic.

Furthermore, we introduce the $k$-jets of $d$-dimensional submanifolds.
Let $\c_x^d(X)$ be the set of all germs at $x$
of $d$-dimensional submanifolds of $X$.
We say that two germs $V,V'\in\c_x^d(X)$
are $k$-equivalent, if, in some local coordinate neighborhoods $U_1\times U_2$ of $x$,
$V$ and $V'$ can be given as graphs of
holomorphic maps $\phi,\phi'\colon U_1\to U_2$
such that $j^k_{x_1}\phi=j^k_{x_1}\phi'$.
Denote by $J_x^{k,d}(X)$ the space of all $k$-equivalence classes at $x$
and by $J^{k,d}(X)$ the union $\cup_{x\in X}J^{k,d}(X)$
with the natural fiber bundle structure over $X$.
Furthermore, for $V\in\c_x^d(X)$, denote by $j^k_x(V)\in J_x^{k,d}(X)$
the corresponding $k$-jet.
For $g\in J^{k,d}(X)$, set $g(0):=x$, if $g$ is a $k$-jet at $x$.
If $X$ is an algebraic manifold, $J_x^{k,d}(X)$ and $J^{k,d}(X)$
are also algebraic.

\begin{Def}\label{inc}
Let $g\in J^{k,d}(X)$ and $g'\in J^{k,d'}(X)$ be arbitrary $k$-jets.
We say that $g$ is contained in $g'$ and write $g\subset g'$,
if $g(0)=g'(0)$ and there exist germs $V\in g$ and $V'\in g'$ such that
$V\subset V'$.
\end{Def}

Let $F\subset U\times V$ and $F'\subset U'\times V$
be holomorphic families. Analogously to $A(F,F')$ we define
$$ A^k(F,F'):= \{ (u,v,u') \in U\times V\times U'
	: j^k_v F'_{u'} \supset j^k_v F_u \}.$$
It follows that $A^k_{u,v}(F,F')\supset A^{k+1}_{u,v}(F,F')$.
Since the fibers $F_u$ are connected,
$$A_u(F,F')=\cap_k A^k_{u,v}(F,F'), \quad (u,v)\in F.$$
Hence $\pi^{-1}(A(F,F'))$ equals to the intersection
of the decreasing sequence of the analytic sets $A^k(F,F')$,
where $\pi\colon U\times V\times U' \to U\times U'$
denotes the natural projection. The following property
is straightforward:

\begin{Lemma}\label{akff}
Let $F\subset U\times V$ and $F'\subset U'\times V$ be holomorphic families,
$(u,v)\in F$ and $(u',v)\in F'$ be such that $F_u\subset F'_{u'}$.
Then we can replace $F$ and $F'$ with subfamilies
containing $(u,v)$ and $(u',v)$ respectively and such that
$$\pi^{-1}(A(F,F'))=A^k(F,F')$$
for $k$ sufficiently large.
\end{Lemma}

Furthermore, for $0\le d\le \dim V$ we define
$$A^{k,d}_g(F):= \{u\in U : j^k_{g(0)} F_u\supset g\},$$
where $g\in J^{k,d}(V)$ is an arbitrary jet.

\begin{Lemma}\label{akdf}
Let $F\subset U\times V$ be a holomorphic (resp. algebraic) family.
Then $A^{k,d}(F)\subset J^{k,d}(V)\times U$ is an analytic
(resp. algebraic) subset.
\end{Lemma}

\begin{proof}
Let $\phi\colon U\times V_1\to V_2$ be as in Definition~\ref{family}.
Suppose that $(g_0,u_0)\in A^{k,d}(F)$ is arbitrary.
Then $d\le \dim V_1$ and there exists a splitting $V_1=W_1\times W_2$
such every $k$-jet $g\in J^{k,d}(V)$ sufficiently closely to $g_0$
can be represented by a graph of an algebraic map
$(\phi_g,\psi_g)\colon W_1\to W_2\times V_2$.
Now the condition  $(g,u)\in A^{k,d}(F)$ for $(g,u)$ close to $(g_0,u_0)$
is equivalent to
$$\partial^{|\alpha|}\phi(u,w_1,\phi_g(w_1))/\partial w_1^{\alpha}
= \partial^{|\alpha|}\psi_g(w_1)/\partial w_1^{\alpha}$$
for all multi-indices $\alpha$ with $|\alpha|\le k$,
which is an analytic (resp. algebraic) condition
on $u$ and the coordinates of $g$.
\end{proof}

\section{Reflection of jets}\label{reflection}

In this section we consider real algebraic submanifolds
$M\subset\C^n$ and $M'\subset\C^{n'}$
and suppose that the conditions of Theorem~\ref{main} are satisfied
and $f\colon D\to\C^{n'}$ is a holomorphism between $M$ and $M'$.
In addition we also suppose that $M$ and $M'$ are generic.
After replacing $D$ with a smaller domain
we may assume that $M$ is given by (\ref{Meq}) with $U=D$, $p_i$ real-valued and
\begin{equation}\label{wed}
(\partial p_1(z,\w)/ \partial z)\wedge\cdots\wedge (\partial p_m(z,\w)/\partial z)\ne 0
\end{equation}
for all $z,w\in D$.
Then the complexification of $M$ can be defined as follows:
$$\M:=\{(\w,z)\in D\times\ov{D} : p_1(z,\w)=\cdots=p_m(z,\w)=0\}.$$
In particular, $\M$ is an algebraic manifold.
Similarly we may assume that the complexification $\M'$ of $M'$ is also
an algebraic manifold. Then the complexification
of the condition $f(M)\subset M'$ yields $(\f,f)(\M)\subset\M'$,
which is equivalent to
\begin{equation}\label{reflect}
f(Q_w)\subset Q'_{f(w)}, \, w\in D,
\end{equation}
where $Q_w:=\{z : (\w,z)\in\M \}$ and $Q'_{w'}:=\{z' : (\w',z')\in\M' \}$
are the Segre varieties associated to $M$ and $M'$ respectively.
In the sequel we identify $M$ (resp. $M'$) with its diagonal embedding
$z\mapsto (\z,z)$ into $\M$ (resp. $\M'$).

Let $X_i$ (resp. $X'_i$) be domains in $\cn$ (resp. $\cnn$) for $i$ odd
and in $\ocn$ (resp. $\ocnn$) for $i$ even.
Set $\M_i:= \M\cap (X_{i-1}\times X_i)$
(resp. $\M_i:= \ov{\M}\cap (X_{i-1}\times X_i)$) 
 for $i$ odd (resp. for $i$ even).
We define the iterated complexifications of $M$ as follows:
\begin{equation}\label{Mm}
\M^m:=\{\xi=(\xi_0,\ldots,\xi_m)\in X_0\times\cdots\times X_m
	: (\xi_0,\xi_1)\in \M_1, \ldots, (\xi_{m-1},\xi_m)\in \M_m \}.
\end{equation}
We also use the notation $f_i(\xi):=f(\xi_i)$
for $i$ odd and $f_i(\xi):=\f(\xi_i)$ for $i$ even.
Of course, this construction depends on the choice of the domains
$X_i$, $i=1,\ldots,m$.

\begin{Def}\label{admis}
We call the tuple $X(m)=(X_1,\ldots,X_m,X'_0,\ldots,X'_m)$
{\em admissible}, if the following is satisfied:
\begin{enumerate}
\item $X_i\subset D$ (resp. $X_i\subset \ov{D}$) for $i$ odd (resp. for $i$ even);
\item $f(X_i)\subset X'_i$ (resp. $\f(X_i)\subset X'_i$) for $i$ odd (resp. for $i$ even);
\item there exists $x\in M$ such that 
	$(\x,x,\x,\ldots)\in X_0\times X_1\times\cdots\times X_m$.
\end{enumerate}
\end{Def}

Clearly we obtain an admissible tuple $X(m)$ if we set $X_i:=D$ and $X'_i:=\cnn$
for $i$ odd (resp. $X_i:=\ov{D}$ and $X'_i:=\ocnn$ for $i$ even).
By a {\em subtuple} of $X(m)$ we mean a tuple
$\tilde X(m)=(\tilde X_1,\ldots,\tilde X_m,\tilde X'_0,\ldots,\tilde X'_m)$
such that $\tilde X_i\subset X_i$ and $\tilde X'_i\subset X'_i$.
In the sequel we shall frequently replace $X(m)$
with subtuples which we also denote by $X(m)$ if no ambiguity appears.

\begin{Lemma}\label{smooth}
Let $U$, $V$ be algebraic manifolds, $S\subset U\times V$ be an algebraic subset
and $M\subset S$ be a real analytic submanifold.
Then there exist an algebraic submanifold $N\subset U$,
an open subset $V'\subset V$ and an algebraic family $F\subset N\times V'$ such that
$F\subset S$ and $F\cap M$ is open in $M$.
\end{Lemma}

\begin{proof}
We prove the statement by induction on $\dim S$.
It is clear for $S$ discrete.
If $M\subset \sing S$, the statement follows by the induction hypothesis,
because $\dim\sing S< \dim S$.
Otherwise we may assume that $S$ is smooth by replacing it with $\reg S$.
By a similar argument we may assume that the natural projection
$\pi_{S,U}\colon S\to U$ is of constant rank.
Then $\pi_{S,U}$ is locally a submersion onto an algebraic submanifold $N\subset U$
and the statement follows from Lemma~\ref{submer}.
\end{proof}

\begin{Prop}\label{psi}
Under the assumptions of Theorem~{\rm\ref{main}}
suppose that $M$ and $M'$ are generic.
Then $X(2)$ can be replaced with an admissible subtuple such that
for every nonempty open set $\Omega\subset\M^2$ and $k$ sufficiently large, 
there exist a nonempty open subset $\Omega'\subset\Omega$ and a partial algebraic map
$\Psi\colon \M^2\times J^k(X_1,X'_1) \times J^k(X_0,X'_0)\to X'_0$ with
\begin{equation}\label{Psi}
\f(\w_1)=\Psi(\xi,j^k_zf,j^k_\w\f)
\end{equation}
for all $\xi=(\w_1,z,\w)\in \Omega'$.
\end{Prop}

\begin{proof}
We apply the properties obtained in previous sections to our situation.
Set $U:=X_0$, $U':=X'_0$, $V:=X'_1$, $F_\w:=f(Q_w)$, $F'_{\w'}:=Q'_{w'}$.
By Lemmata~\ref{submer} and ~\ref{image} we can replace $X(2)$
with an admissible subtuple such that
$F\subset U\times V$ and $F'\subset U'\times V$
become holomorphic families. Since $\M'$ is algebraic, $F'$ is an algebraic family.
For $d:=\dim f(Q_w)$ define:
$$ A^k := A^{k,d}(F')\subset J^{k,d}(X'_1)\times X'_0.$$
By Lemma~\ref{akdf}, $A^k$ is an algebraic subset.
By Lemma~\ref{akff}, $X(2)$ can be replaced with a subtuple such that for all $(\w,z)\in\M$
and $k$ sufficiently large,
\begin{equation}\label{ak}
A^k_{g(\w,z)}=A_\w(F,F') =
	\{ \w'\in X'_1 : Q'_{w'} \supset f(Q_w) \},
\end{equation}
where $g(\w,z):=j^k_{f(z)}f(Q_w)\in J^{k,d}(X'_1)$.
Together with (\ref{reflect}) this implies
\begin{equation}\label{a-inc}
\f(\w)\in A^k_{g(\w,z)}, \quad (\w,z)\in\M.
\end{equation}

Set $\tilde \M:=\{(g(\w,z),\f(\w)): (\w,z)\in\M\}$
and $\tilde M:=\{(g(\z,z),\f(\z)): z\in M\}$.
By Lemma~\ref{smooth}, we can replace $A^k$ with an algebraic subset $A\subset S$
which is a family over 
some algebraic submanifold $N\subset J^{k,d}(X'_1)$ such
that $\tilde M\cap A$ is open in $\tilde M$. Since $\tilde M$ is generic in $\tilde\M$
and $A$ is complex, $\tilde\M\cap A$ is open in $\tilde\M$.
Hence we can replace $X(2)$ with an admissible subtuple such that (\ref{a-inc}) is still valid.

Now consider the conjugated family $\A\subset J^{k,d}(\ov{X'_1})\times \ov{X'_0}$.
By replacing $X(2)$ with an admissible subtuple we may assume that $\A$ is an algebraic
family in $\ov N\times X'_1$. By Lemma~\ref{bag},
$$ B := B(\A,\ov{F'})\subset J^{k,d}(X'_0)\times X'_0$$
is an algebraic subset. Since $A\subset A^k$, it follows from (\ref{ak}) that
$$\f(\ov{Q_w}) \subset \cap_{w'\in \A_{\g(w,\z)}} \ov{Q'_{w'}} = B_{\g(w,\z)}, \quad  (\w,z)\in\M.$$
By the construction of $\M$, $(\w,z)\in\M$ is equivalent to $(\z,w)\in\M$
for $z,w$ close to $M$.
Hence, by replacing $X(2)$ with an admissible subtuple we may assume that
\begin{equation}\label{b-inc}
\f(\ov{Q_z}) \subset B_{\g(z,\w)}, \quad  (\w,z)\in\M.
\end{equation}
Furthermore, we set $J_1:=J^{k,d}(X'_1)$ and $J_2:=J^{k,d}(X'_0)$ and define:
$$C:=\{(g_1,g_2,\w')\in J_1\times J_2\times X'_0 : 
	\w'\in A_{g_1}\cap B_{g_2} \}.$$
Let $\xi=(\w,z,\w_1)\in\M^2$ be arbitrary.
It follows from (\ref{a-inc}) and (\ref{b-inc}) that
\begin{equation}
\f(\w_1)\in C_{g(\w_1,z),\g(z,\w)}.
\end{equation}

Since $X(2)$ is admissible, there exists $x\in M$ as in Definition~\ref{admis},
i.e. $(\x,x,\x)\in\M^2$. Then 
$$ C_x:= C_{g(\x,x),\g(x,\x)} \subset B_{\g(x,\x)} \subset \ov{Q'_{w'}}$$
for all $\w'\in C_x\subset \A_{\g(x,\x)}$.
In particular, this implies $\w'\in \ov{Q'_{w'}}$, i.e. $(\w',w')\in\M'$.
By the construction of $\M'$, this is equivalent to $\w'\in M'$.
Thus the whole analytic fiber $C_x$ lies in $M'$.
By Condition~2 of Theorem~\ref{main}, we conclude that $C_x$ must be discrete.
Since $C\subset J_1\times J_2\times X'_0$ is an algebraic subset,
all fibers $C_{g_1,g_2}$ in a neighborhood of $(g(\x,x),\g(x,\x),\f(\x))\in C$ 
must be discrete. 

Let $\Omega\subset\M^2$ be an open subset.
By taking $\Omega$ smaller we may assume that
the map $(\w_1,z,\w)\mapsto (g(\w_1,z),\g(z,\w),\f(w_1))$
sends $\Omega$ onto a complex submanifold 
$\tilde\Omega\subset C$.
By Lemma~\ref{smooth}, there exists a nonempty open subset 
$\Omega'\subset\Omega$ and a partial algebraic map
$\Phi\colon J_1\times J_2\to X'_0$ such that
\begin{equation}\label{Phi}
\f(\w_1) = \Phi(g(\w_1,z),\g(z,\w)), \quad (\w_1,z,\w)\in\Omega'.
\end{equation}

Set $e:=\dim Q_w$, $\w\in X_0$. 
The jet composition yields a partial algebraic map
$$\phi^k\colon J^k(X_1,X'_1) \times J^{k,e}(X_1) \to J_1=J^{k,d}(X'_1)$$
such that $j^k_{f(z)}f(Q_w)=\phi^k(j^k_zf, j^k_z Q_w)$
for all $(\w,z)\in\M$. Together with (\ref{Phi}) we obtain:
\begin{equation}
\f(\w_1)=\Phi\left(\phi^k(j^k_zf,j^k_z Q_{w_1}),
	\ov{\phi^k}(j^k_\w\f,j^k_\w \ov{Q_z})\right).
\end{equation}
This proves the required statement.
\end{proof}

\section{Iteration of jet reflections}

\begin{Prop}\label{psi-jet}
Under the assumptions of Theorem~{\rm\ref{main}}
suppose that $M$ and $M'$ are generic.
Then we can replace $X(2)$ with an admissible subtuple such that
for every nonempty open set $\Omega\subset\M^2$ and $k$ sufficiently large, 
there exists a nonempty open subset $\Omega'\subset\Omega$ 
and partial algebraic maps
$$\Psi^s\colon \M^2\times J^{s+k}(\M^2,X'_1) \times
J^{s+k}(\M^2,X'_0)\to
    J^s(\M^2,X'_2), \quad s=0,1,\ldots,$$
such that $j^s_\xi f_2=\Psi^s(\xi,j^{s+k}_\xi f_1,j^{s+k}_\xi f_0)$ for all
$\xi=(\w_1,z,\w)\in \Omega'$.
\end{Prop}

\begin{proof}
The statement is obtained by taking the $s$-jets of both sides in (\ref{Psi}).
\end{proof}

\begin{Prop}\label{Ml}
Under the assumptions of Theorem~{\rm\ref{main}}
suppose that $M$ and $M'$ are generic.
Then for every $l=1,2,\ldots$ we can replace $X(l)$ with an admissible subtuple such that
for every nonempty open set $\Omega\subset\M^l$ and $k$ sufficiently large, 
there exists a nonempty open subset $\Omega'\subset\Omega$ 
and partial algebraic maps
$$H^s_l\colon \M^l\times J^m(\M^l,X'_1) \times J^m(\M^l,X'_0)
	\to J^s(\M^l,X'_l),$$
such that
\begin{equation}\label{ml}
j^s_\xi f_l=H^s_l(\xi,j^m_\xi f_1,j^m_\xi f_0)
\end{equation}
for all $\xi\in \Omega'$, $s\ge 0$, $l\ge 1$ and $m:=s+k(l-1)$.
\end{Prop}

\begin{proof}
We argue by induction on $l$.
For $l=1$ the statement is trivial, for $l=2$
it is contained in Proposition~\ref{psi-jet}.
We now suppose that the statement is true for $l-2$ and $l-1$.
Define
$$\M^2_l:=\{(\xi_l,\xi_{l-1},\xi_{l-2})\in X_l\times X_{l-1}\times X_{l-2} :
	(\xi_l,\xi_{l-1})\in\M_l, (\xi_{l-1},\xi_{l-2})\in\M_{l-1} \}.$$
Then $\M^2_l$ is either canonically biholomorphic or anti-biholomorphic
to $\M^2$. By Proposition~\ref{psi-jet}, we can replace 
$(X_{l-2},X_{l-1},X_l,X'_{l-2},X'_{l-1},X'_l)$ with an admissible subtuple
such that for every nonempty open set $\Omega_l\subset\M^2_l$ and $k$ sufficiently large,
there exists a nonempty open subset $\Omega'_l\subset\Omega_l$ and
partially algebraic maps
$$P^s_l\colon \M^2_l\times J^{s+k}(\M^2_l,X'_{l-1})
	\times J^{s+k}(\M^2_l,X'_{l-2})\to J^s(\M^2_l,X'_l), \quad
       s=0,1,\ldots,$$
such that
\begin{equation}\label{chi}
j^s_\chi f'_l=P^s_l(\chi,j^{s+k}_\chi f'_{l-1}, j^{s+k}_\chi f'_{l-2}),
\quad \chi=(\xi_{l-2},\xi_{l-1},\xi_l)\in\Omega'_l,
\end{equation}
where $f'_i(\chi):=f(\xi_i)$ for $i$ odd and 
$f'_i(\chi):=\f(\xi_i)$ for $i$ even.
Then $f_i(\xi)=f'_i(\xi_{l-2},\xi_{l-1},\xi_l)$ and the above property
can be reformulated as follows:
{\em For every nonempty open set $\Omega\subset\M^l$
there exists a nonempty open subset $\Omega'\subset\Omega$
such that
\begin{equation}\label{chi'}
j^s_\xi f_l=\tilde P^s_l(\xi,j^{s+k}_\xi f_{l-1}, j^{s+k}_\xi f_{l-2}),
\quad \xi\in\Omega',
\end{equation}
where 
$$\tilde P^s_l\colon \M^l\times J^{s+k}(\M^l,X'_{l-1})
	\times J^{s+k}(\M^l,X'_{l-2})\to J^s(\M^l,X'_l),$$
are partial algebraic maps}.

On the other hand, we can apply the induction hypothesis for $(l-1)$ and \hbox{$(l-2)$}.
Without loss of generality, $X(l)$ satisfies the properties given by
the inductive assumptions. This means that we can replace 
the above set $\Omega'$ with a nonempty open subset such that
\begin{equation}\label{i}
j^{s+k}_\xi f_i=\tilde H^{s+k}_i(\xi,j^{m+k}_\xi f_1,j^{m+k}_\xi f_0),
\quad \xi\in\Omega',
\end{equation}
where
$$\tilde H^s_i\colon \M^l\times J^m(\M^l,X'_1) \times J^m(\M^l,X'_0)
	\to J^s(\M^l,X'_i)$$
are some partial algebraic maps for $i=(l-1),(l-2)$.
The statement follows from (\ref{chi'}) and (\ref{i}).
\end{proof}

\begin{proof}[Proof of Theorem~{\rm\ref{main}}]
We first assume that $M$ and $M'$ are generic.
Using the projection $\pi_i\colon \M^l\to X_i$ we obtain the canonical jet lifting 
$$J^s(X_i,X'_i)\to J^s(\M^l,X'_i),\quad g\mapsto \pi_i^*g, \quad i=0,1,\ldots,l.$$
Then we can rewrite (\ref{ml}) for $s=0$ and $l$ odd as follows:
\begin{equation}\label{fl}
f(\xi_l)=H(\xi, j^m_{\xi_1} f, j^m_{\xi_0} \f), \quad \xi\in\M^l,
\end{equation}
where
$$H\colon \M^l\times J^m(X_1,X'_1) \times J^m(X_0,X'_0)	\to X'_l$$
is a partial algebraic map.

Now we use the following consequence of Theorem~3 and Corollary~2.2.2 in \cite{BER}:

\begin{Cor}\label{gen-sub}
Under the assumptions of Theorem~\ref{main}, 
$\pi_l\colon \M^l_{\chi_1,\chi_0} \to X_l$
is generically submersive for $l$ sufficiently large, where
$\M^l_{\chi_1,\chi_0}:= \{\xi\in\M^l : \xi_1=\chi_1, \xi_2=\chi_2 \}$.
\end{Cor}

Without loss of generality, $l$ is odd.
Then there exists an open subset $U\subset X_l$ and an algebraic
lifting $\psi\colon U\to \M^l_{\chi_1,\chi_0}$ such that $\pi_l\circ\psi=\id_U$.
Therefore, (\ref{fl}) implies
\begin{equation}
f(\xi_l)=H(\psi(\xi_l), j^m_{\chi_1} f, j^m_{\chi_0} \f), \quad \xi_l\in U.
\end{equation}
Hence $f$ is algebraic as composition of the algebraic maps $\psi$ and $H$.
\end{proof}

It remains to reduce the general case to the case of generic $M$ and $M'$.
Suppose that $f$ is defined in $D\subset\cn$, $x\in M\cap D$ and $x':=f(x)\in M'$.
Recall that $M\subset\cn$ is said to be a \CR submanifold,
if $\dim T_xM\cap iT_xM$ does not depend on $x\in M$.
We apply the following lemma (the proof is analogous to that of Lemma~\ref{smooth}):

\begin{Lemma}\label{smooth1}
Let $M\subset\cn$ be a real-algebraic subset and let $S\subset M$ be a real-analytic submanifold.
Then there exists an algebraic \CR submanifold $N\subset M$,
such that $N\cap S$ is open in $S$.
\end{Lemma}

Without loss of generality, $M$ and $M'$ are connected
and $f|M$ is of constant rank.
By applying Lemma~\ref{smooth1} to $M\cap D\subset M$
and $f(M\cap D)\subset M'$ we may assume that $M$ and $M'$ are \CR submanifolds.
Then there exist algebraic local coordinates $(z_1,z_2)\in\cn$ near $x$ and 
$(z'_1,z'_2)\in\cnn$ near $x'$ such that $M$ (resp. $M'$) is generic in 
$\{z_2=0\}$ (resp. $\{z'_2=0\}$). However, condition~1 of Theorem~\ref{main}
implies that $z_2=z$ (otherwise $(z_1,z_2)\mapsto \Re\, z_2$ would be a nonconstant 
holomorphism between $M$ and $\R$). Hence $M$ is already generic.
Therefore the condition $f(M\cap D)\subset M'\subset \{z'_2=0\}$ implies $f(D)\subset\{z'_2=0\}$.
Thus we can apply the statement to $M\subset\cn$ and $M'\subset \{z'_2=0\}\equiv \C^{n''}$.
Since the local coordinates $(z'_1,z'_2)$ are algebraic, we obtain the required
algebraicity of $f$.

\section{An algorithm for finding analytic discs in a real submanifold}\label{algo}

Let $M\subset\cn$ be a real submanifold given locally by the equation
\begin{equation}\label{meq}
M= \{ z\in U : r(z,\z)=0 \},
\end{equation}
where $r\colon U\times\bar U\to \R^k$ is a real-valued convergent power series.
Let $x\in M$ be an arbitrary point. By an analytic disc in $M$ we mean
a nonconstant holomorphic map $f$ from the unit disc $\Delta:=\{|z|<1\}\subset\C$ into $\C^n$ such
that $f(\Delta)\subset M$.

{\bf Step 1.}
Define the complexification as above:
\begin{equation}\label{c-tion}
\M:=\{ (z,\w)\in U\times\bar U : r(z,\w)=0 \}.
\end{equation}
Moreover, we assume that $r$ is analytic on the closure of $U\times\bar U$.
This will guarantee that the analytic subsets defined in terms of $r$ 
have finitely many irreducible components in $U\times\bar U$.

{\bf Step 2.}
Construct the family of intersections:
\begin{equation}
A:=\{ (a_1,\ldots,a_k,\w)\in U^k\times\bar U : r(a_1,\w)=\cdots=r(a_k,\w)=0 \},
\end{equation}
where $k:=[\dim_\R M/2]$. Define $A_a:=\{\w\in\bar U: (a,\w)\in A\}$.

{\bf Step 3.}
Determine all points $z$ such that $r(z,\cdot)$ vanishes on $A_a$:
\begin{equation}\label{C2}
B:=\{ (a,z)\in U^k\times U : r(z,A_a)=0 \}.
\end{equation}

{\bf Step 4.}
Define
\begin{equation}\label{C1}
C:=\{ (a,z)\in B : (a,\z)\in A \}, \quad C_a:=\{z\in U: (a,z)\in C\}.
\end{equation}

Then $M\cap U$ contains no nonconstant complex-analytic discs 
if and only if $\dim C_a=0$ for all $a\in U^k$ according to the following proposition: 

\begin{Prop}
For every $a\in U^k$, $C_a$ is a complex-analytic set which is contained in $M\cap U$.
On the other hand, every complex-analytic disc in $M\cap U$ is contained in $C_a$ for some $a\in U^k$.
\end{Prop}

\begin{proof}
The analyticity of $C_a$ is straightforward. 
If $z\in C_a$ is arbitrary, $\z\in A_a$ by (\ref{C1}) and therefore $r(z,\z)=0$ by (\ref{C2}).
It follows from (\ref{meq}) that $z\in M$. This proves that $C_a\subset M$.

On the other hand, let $f\colon \Delta\to M\cap U$ be a complex-analytic disc.
The reflection property (\ref{reflection}) yields 
$f(\Delta)\subset Q_{f(t)}$ for all $t\in\Delta$,
where $Q_w:=\{z\in U: r(z,\w)\}$ is the Segre variety attached to $M$.
Hence $f(\Delta)$ is contained in the intersection $I$ of all $Q_{f(t)}$, $t\in\Delta$.
The set $I$ is of complex dimension not larger than $k=[\dim_\R M/2]\ge \dim_\C Q_{f(t)}$
and has finitely many irreducible components by the assumption in Step~1.
Thus $I=A_a$ for some $a\in f(\Delta)^k$. 
Since $A_a=I\subset Q_{f(t)}$, we have $(a,f(t))\in B$ for all $t\in\Delta$. 
On the other hand, $f(t)\in A_a$ by our construction.
This proves that $f(\Delta)\subset C_a$.
\end{proof}


\end{document}